\documentclass[a4paper] {article} [12pt]
\usepackage{amsmath,amsfonts,amssymb,amsthm,amscd,nccmath}
\usepackage{graphicx}
\usepackage[numbers]{natbib}
\usepackage{color}
\usepackage{xypic}
\usepackage{tikz}
\usepackage{environ}
\usepackage{fullpage}

\newtheorem{Theorem}{Theorem}[section]
\newtheorem{thm}[Theorem]{Theorem}

\newtheorem{defn}[Theorem]{Definition}

\newtheorem{rem}[Theorem]{Remark}

\newtheorem{lem}[Theorem]{Lemma}
\newtheorem{Lem}[Theorem]{Lemma}

\newtheorem{cor}[Theorem]{Corollary}

\newtheorem{fact}[Theorem]{Fact}

\makeatletter
\newsavebox{\measure@tikzpicture}
\NewEnviron{scaletikzpicturetowidth}[1]{%
 \def\tikz@width{#1}%
 \begin{lrbox}{\measure@tikzpicture}%
 \BODY
  \end{lrbox}%
  \pgfmathparse{#1/\wd\measure@tikzpicture}%
  \BODY
}
\makeatother

\makeatletter
\newcommand{\thorn}{{\fontencoding{T1}\selectfont\th}}
\newcommand{\@indepsymbol}[2]{#1\setbox0=\hbox{$#1x$}\kern\wd0\hbox to 
0pt{\hss$#1\mid$\hss}\lower.9\ht0\hbox to 0pt{\hss$#1\smile$\hss}\kern\wd0}
\newcommand{\@nindepsymbol}[2]{#1\setbox0=\hbox{$#1x$}\kern\wd0\hbox to 
0pt{\mathchardef
	\nn=12854\hss$#1\nn$\kern1.4\wd0\hss}\hbox to
	0pt{\hss$#1\mid$\hss}\lower.9\ht0 \hbox to
	0pt{\hss$#1\smile$\hss}\kern\wd0}
\newcommand{\ind}[1][]{\mathop{\mathpalette\@indepsymbol{}^{\!\!\!\!\rlap{
$\scriptstyle\textnormal{#1}$}\,\,\,\,}}}
\newcommand{\nind}[1][]{\mathop{\mathpalette\@nindepsymbol{}^{\!\!\!\rlap{
$\scriptstyle\textnormal{#1}$}\,\,\,}}}
\newcommand{\@Ind}[1][]{\mathpalette\@indepsymbol{}^{\!\!\!\!\mbox{
$\scriptstyle\textnormal{#1}$}}}
\newcommand{\Ind}[1][]{\@Ind[\ \,]}
\newcommand{\newind}[4]{
	\newcommand{#1}{{\!\@Ind[#4]}}
	\newcommand{#2}{\ind[#4]}
	\newcommand{#3}{\nind[#4]}
}
\newind{\thInd}{\thind}{\nthind}{\thorn}
\renewcommand{\thInd}{\text{$\@Ind[\thorn]$\;}}
\makeatother

\title{Theories with distal Shelah expansions}

\date{\today}
\author{Gareth Boxall \thanks{This work is based on the research supported in part by the National Research Foundation of South Africa (Grant number 96234).} \and Charlotte Kestner\thanks{Supported by LMS travel grant 41605 and Leverhulme Project grant RPG-2017-179.}}

\begin{document}

\maketitle

\begin{abstract}
We show that a complete first-order theory $T$ is distal provided it has a model 
$M$ such that the theory of the Shelah expansion of $M$ is distal.  
\end{abstract}

\section{Introduction}

Since its introduction \cite{Sdist} distality has played an important 
role in the study of NIP. The notion identifies those NIP theories 
which are, in some sense, completely non-stable. Thus RCF is distal while ACVF, which has 
o-minimal value group and stable residue field, is not. Indeed, Simon has 
decomposition results \cite{NIPtypedecomp} according to which types in NIP theories are 
understood in terms of a stable part and a distal part. The theory of an 
infinite set is not distal and so distality has the quirky property, among 
tameness conditions, of not always passing to reducts. 

Of course, some reducts of structures with distal theories 
will have distal theories and so will some expansions. In the NIP context, it is natural to 
consider the Shelah expansion: one adds to the language a predicate for every 
externally definable set of the structure. Shelah proved \cite{ShExp} that NIP is preserved 
when moving to this expansion and trivially NIP passes to reducts. In early 
2017, Artem Chernikov pointed out to us that, while it is easy to see that the 
Shelah expansion of a  model of a distal theory will have a distal theory, it is 
not so clear that the intermediate expansions will have distal theories. We prove that they will. 

\begin{thm}\label{main}
Let $T$ be a complete first-order theory. Let $M\models T$ and let $M^{Sh}$ be 
the Shelah expansion of $M$. If $Th(M^{Sh})$ is distal then $T$ is distal.
\end{thm}

\begin{cor}\label{question}
If $\hat{M}$ is an expansion of $M$ such that $M^{Sh}$ is an expansion of $\hat{M}$, then distality of $Th(M^{Sh})$ implies distality of $Th(\hat{M})$. 
\end{cor}

In the next section we say exactly what we mean by ``expansion'' and show how the corollary is obtained from the theorem. Our proof of Theorem \ref{main} relies on a lemma in the NIP setting which we hope will have other applications. Before stating it, we mention some notational conventions. We identify a non-constant indiscernible sequence $(a_i)_{i\in I}$ with the ordered set $\{a_i:i\in I\}$, the order given by saying $a_i<a_j$ if and only if $i<j$. If $M$ is a structure and $A$ a set then $A\subseteq M$ means $A$ is a subset of a finite cartesian power of sorts of $M$. If we want to be specific about which cartesian power then we write $A\subseteq M^{\bar{x}}$, where $\bar{x}$ is a tuple of variables of the appropriate sorts. We follow a similar convention with the notation $a\in M$ and $a\in M^{\bar{x}}$. If $A$ and $B$ are both just sets then $A\subseteq B$ has its usual meaning. 

We make much use of pairs of structures $(N,M)$ where $M\prec N$. In the one-sorted setting, the language of such a structure is $L_P=L\cup\{P\}$, where $L$ is the language of $M$ and $P$ is a new unary predicate interpreted such that $P(N)=M$. In the many-sorted setting, one would need to replace $P$ with a family $(P_s)_{s\in S}$ of unary predicates, one for each sort, and the interpretation would be $P_s(N_s)=M_s$ for each $s\in S$. Given that this is understood, we shall for simplicity use the one-sorted notation even in the many sorted setting. 

\begin{Lem}\label{extending}
Let $L$ be a language and $T$ a complete first-order $L$-theory. Assume $T$ has 
NIP. Let $M\models T$ and let $M\prec N$ and $(N,M)\prec (N',M')$ be 
sufficiently saturated elementary extensions, where $(N,M)$ and $(N',M')$ are $L\cup\{P\}$-structures and $P$ is a unary predicate, not in $L$, such that $P(N)=M$ and $P(N')=M'$. Let $A\subseteq M'$ be a small 
non-constant $L$-indiscernible sequence. Then there is a small 
$L$-indiscernible sequence $A'\subseteq M'$ which extends $A$ and has the following property. For every complete 
$L$-type $q(\bar x)$ over $A'$ with $\bar{x}=(x_1,...,x_n)$, if $q(\bar{x})$ is finitely realised in $A'$ then $q(\bar{x})\cup \{P(x_1),...,P(x_n)\}$ implies a complete $L$-type over $N$.
\end{Lem}

This lemma is one of two ingredients in the proof of Theorem \ref{main}. The other is the argument used in \cite{CS} to show that distality is equivalent to the existence of strong honest definitions (Theorem 21 in \cite{CS}). 

In \S 2 we give some background definitions 
and information. In \S 3 we prove Lemma \ref{extending}. In \S 4 we prove 
Theorem \ref{main}. In \S 5, for completeness, we prove the converse of 
Theorem \ref{main} which was already known to experts. 

After talks given by the second author on this work, Ehud Hrushovski and Anand Pillay both directed us to an alternative approach to Theorem \ref{main} via generically stable measures, using Simon's characterisation of distal theories as those NIP theories for which every generically stable measure over a model is smooth. Both were kind enough to supply us with further details and to grant permission for the inclusion of the argument here. We sketch it in \S 6. We felt it appropriate to retain our original proof, which avoids use of measures and Simon's result, partly as an advertisement for Lemma \ref{extending}. 

We thank Artem Chernikov for alerting us to the problem addressed in this paper and for helpful 
discussions. We thank Ehud Hrushovski and Anand Pillay for \S6. We thank the referee for a thorough reading and several comments and suggestions which helped correct errors and improve the paper.

\section{Preliminaries}\label{prelim}

Let $L$ be a language and $T$ a complete first-order $L$-theory. Let $M\models T$ and let $M\prec M'$ be sufficiently saturated. We assume throughout that $T$ has NIP. By an expansion of $M$ we mean a structure $N$ with the same underlying set (the same sorts, in the many sorted-setting) such that every $\emptyset$-definable set of $M$ is $\emptyset$-definable in $N$. An externally definable set of $M$ is a set of the form $X\cap M^{\bar{x}}$ where $X\subseteq M'^{\bar{x}}$ is definable (with parameters) in $M'$. It is easy to check that this does not depend on the choice of saturated extension $M'$. The Shelah expansion $M^{Sh}$ of $M$ is the structure whose language $L(M^{Sh})$ has one predicate for each externally definable set of $M$ and in which these predicates are interpreted in the obvious way. We shall rely on the following fact proved by Shelah in \cite{ShExp}. 

\begin{fact}\label{QE}
$Th(M^{Sh})$ has quantifier elimination and NIP. 
\end{fact}

Let $\hat{M}$ be an expansion of $M$ such that $M^{Sh}$ is an expansion of $\hat{M}$. Corollary \ref{question} follows from Theorem \ref{main} in combination with the fact that $M^{Sh}$ is an expansion of $\hat{M}^{Sh}$, it being obvious that $\hat{M}^{Sh}$ is an expansion of $M^{Sh}$. The fact that $M^{Sh}$ expands $\hat{M}^{Sh}$ follows from the fact that $M^{Sh}$ expands $(M^{Sh})^{Sh}$. This must be well known but, unaware of a suitable reference, we provide a short proof (the main points of which were suggested to us by the referee). 

\begin{lem}
$M^{Sh}$ is an expansion of $(M^{Sh})^{Sh}$.
\end{lem}

\proof Let $M^{Sh}\prec \bar{N}$ and $(\bar{N},M^{Sh})\prec (\bar{N}',\bar{M}')$ be sufficiently saturated elementary extensions in the languages $L(M^{Sh})$ and $L(M^{Sh})\cup \{P\}$ respectively. Let $N,M'$ and $N'$ be the reducts of $\bar{N},\bar{M}'$ and $\bar{N}'$ to $L$. Let $X\subseteq M^{\bar{x}}$ be $\emptyset$-definable in $(M^{Sh})^{Sh}$. Then there exists $Y\subseteq M'^{\bar{x}}$ definable in $\bar{M}'$ such that $X=Y\cap M^{\bar{x}}$. We then have some $Z\subseteq M'^{\bar{x}\bar{y}}$ and $\bar{b}\in M'^{\bar{y}}$ such that $Z$ is $\emptyset$-definable in $\bar{M}'$ and $Y=\{\bar{a}\in M'^{\bar{x}}:(\bar{a},\bar{b})\in Z\}$. By Fact \ref{QE}, $Z$ is defined by an $L(M^{Sh})$-formula which defines an externally definable set of $M$. It follows that there exists $W\subseteq N'^{\bar{x}\bar{y}}$ definable in $N'$ such that $Z=W\cap M'^{\bar{x}\bar{y}}$. We then have $X=\{\bar{a}\in M^{\bar{x}}:(\bar{a},\bar{b})\in W\}$ and so $X$ is $\emptyset$-definable in $M^{Sh}$. \endproof

We shall make much use of cuts in the following sense. 

\begin{defn}
Let $(A,<)$ be a totally ordered set. A cut in $A$ is a complete quantifier-free one-type over $A$, considered as a structure in the language $\{<\}$. An unrealised cut in $A$ is one which has no realisation in $A$. 
\end{defn}

The following concept will be useful in the proof of Lemma \ref{extending}. 

\begin{defn}
Let $C=\{c_1,...,c_k\}$ be a set of unrealised cuts in an indiscernible sequence 
$A\subseteq M'$. Let $\bar{b}=(b_1,...,b_n)$ and $\bar{b}^\prime=(b_1^\prime,...,b_n^\prime)$ be 
tuples of elements of $A$. We say that $\bar{b}$ and $\bar{b}^\prime$ have the 
same order type over $C$ if they have the same quantifier-free type over $\emptyset$ in the structure $(A,<)$ and, for each $i\leq n$ and $m\leq k$, $b_i<c_m$ if 
and only if $b_i^\prime<c_m$. In this case we write 
$otp(\bar{b}/C)=otp(\bar{b}^\prime/C)$. 
\end{defn}

The following is an immediate consequence of Fact 1 in \cite{CS} (see also \cite{ShelahMB}). 

\begin{fact} \label{CScutfact} Let $A\subseteq M'$ be a small indiscernible sequence. Let $\theta 
(\bar{x})$ be a formula with parameters in $M'$. Then there is a finite set 
$C=\{c_1,...,c_k\}$ of non-realised cuts in $A$ such that, for any tuples 
$\bar{b},\bar{b}^\prime$ from $A$, if 
$otp(\bar{b}/C)=otp(\bar{b}^\prime/C)$ then $M'\models 
\theta(\bar{b})\leftrightarrow\theta(\bar{b}^\prime)$. 

\end{fact}

We note that there is a minimum such $C$. 

\begin{lem}\label{minimum} Let $A\subseteq M'$ be a small indiscernible sequence. Let $\theta 
(\bar{x})$ be a formula with parameters in $M'$. Let $\mathcal{C}$ be the collection of all $C$ as in Fact \ref{CScutfact}. Then 
$\mathcal{C}$ has a minimum element with respect to set inclusion. 
\end{lem}

\proof By Fact \ref{CScutfact}, $\mathcal{C}$ is not empty. Let 
$C=\{c_1,...,c_k\}\in\mathcal{C}$ such that $k$ is minimal. Let 
$C^\prime\in\mathcal{C}$. Suppose $C\nsubseteq C'$. Let $c\in C\setminus C'$ and $C_1=C\setminus\{c\}$. Since $k$ is minimal, there are tuples $\bar{b},\bar{b}'$ from $A$ such that $otp(\bar{b}/C_1)=otp(\bar{b}'/C_1)$ and $M'\models \theta(\bar{b})\wedge\neg\theta(\bar{b}')$. One can deform $\bar{b}$ into $\bar{b}'$ without changing the truth value of $\theta(\bar{b})$ (by ensuring that at each stage one preserves the order type over $C$ or the order type over $C'$). This is a contradiction. \endproof

Note that none of the cuts in the minimum $C$ will be $\infty$ or $-\infty$. If $A$ and $B$ are disjoint ordered sets then $A+B$ denotes $A\cup B$ equipped with the ordering which places everything in $A$ below everything in $B$ and agrees with the existing 
orderings of $A$ and $B$. The following, which we take as our definition of distality, is provided by a combination of Definition 2.1 and Lemma 2.7 in \cite{Sdist}.  

\begin{defn}\label{distaldef}
$T$ is distal if, for any small indiscernible sequence of the form $I+\{b\}+J$ 
in $M'$, where $\{b\}$ is a singleton and $I$ and $J$ are infinite without endpoints, and any small 
$D\subseteq M'$, if $I+J$ is indiscernible over $D$ then $I + \{b\} + J$ is 
indiscernible over $D$.
\end{defn}

Trivially, if distality fails then this is witnessed by some $D,b,I$ and $J$ such that $I$ and $J$ are both indexed by $\mathbb{Z}$. Note that $M'$ could be many-sorted. Even if it is one-sorted, the elements of $I+\{b\}+J$ could be tuples. When distality fails we shall want the following convenient consequence which must be well known. 

\begin{lem}\label{tech}
If $T$ is not distal then there exist a small indiscernible sequence $I+\{b\}+J$ 
with $I$ and $J$ infinite, a formula $\phi(x,y)$ and some $a\in M'$ such that 
$I+J$ is indiscernible over $a$ and $M'\models \phi(a,c)$ for all $c\in I+J$ but 
$M'\models\neg\phi(a,b)$. 
\end{lem}

\proof Suppose $T$ is not distal. Then we have a small indiscernible 
$I+\{b\}+J$ in $M'$, with $I$ and $J$ indexed by $\mathbb{Z}$, and a small $D\subseteq M'$ such that 
$I+J$ is indiscernible over $D$ but $I+\{b\}+J$ is not. It follows that there 
exist a formula $\phi(x,\bar{y})$ and some $a\in M'$, with $\bar{y}$ an $n$-tuple 
of variables in the sort of $I+\{b\}+J$, such that $M'\models \phi(a,\bar{c})$ 
for any strictly increasing $n$-tuple $\bar{c}$ from $I+J$ and $M'\models 
\neg\phi(a,\bar{b})$ for some strictly increasing $n$-tuple $\bar{b}$ from 
$I+\{b\}+J$. 

Let $I'$ be the set of all elements in $I$ below $\bar{b}$ and let $J'$ be the 
set of all elements in $J$ above $\bar{b}$. Then $I'$ and $J'$ are both infinite and each is indexed by $\mathbb{N}$, with the standard or reverse order, or by $\mathbb{Z}$. By grouping elements together we may treat $I'$ and $J'$ as a sequences of $n$-tuples. Then $I'+\{\bar{b}\}+J'$ is indiscernible and 
$I'+J'$ is indiscernible over $a$. We have $M'\models\phi(a,\bar{c})$, for all 
$\bar{c}\in I'+J'$, while $M'\models\neg\phi(a,\bar{b})$.  \endproof

\section{Lemma}\label{lemmasect}

In this section we prove Lemma \ref{extending}. For convenience we recall the statement.

\begin{Lem}
Let $L$ be a language and $T$ a complete first-order $L$-theory. Assume $T$ has 
NIP. Let $M\models T$ and let $M\prec N$ and $(N,M)\prec (N',M')$ be 
sufficiently saturated elementary extensions, where $(N,M)$ and $(N',M')$ are $L\cup\{P\}$-structures and $P$ is a unary predicate, not in $L$, such that $P(N)=M$ and $P(N')=M'$. Let $A\subseteq M'$ be a small 
non-constant $L$-indiscernible sequence. Then there is a small 
$L$-indiscernible sequence $A'\subseteq M'$ which extends $A$ and has the following property. For every complete 
$L$-type $q(\bar x)$ over $A'$ with $\bar{x}=(x_1,...,x_n)$, if $q(\bar{x})$ is finitely realised in $A'$ then $q(\bar{x})\cup \{P(x_1),...,P(x_n)\}$ implies a complete $L$-type over $N$.
\end{Lem}

\proof When a language other than $L$ is intended, we shall make that clear. Let $\theta(\bar{x})$ be a formula with parameters in $N$. By Fact 
\ref{CScutfact} there is a set $C_A^\theta=\{c_1,...,c_k\}$ of non-realised cuts 
of $A$ such that, for any two tuples $\bar{b}$ and $\bar{b}^\prime$ from $A$, if 
$otp(\bar{b}/C_A^\theta)=otp(\bar{b}^\prime/C_A^\theta)$ then $N'\models 
\theta(\bar{b})\leftrightarrow\theta(\bar{b}^\prime)$. By Lemma \ref{minimum}, 
we may assume $C_A^\theta$ is the minimum among all possible choices (ordered by 
set inclusion). 

For any small indiscernible sequence $B\subseteq M'$ extending $A$ and any cut $c\in 
C_A^\theta$, there must be some $c'\in C_B^\theta$ such that $c^\prime$ refines 
$c$ (by which we mean that every realisation of $c^\prime$ is also a realisation of $c$). Otherwise 
$C_B^\theta$ would give rise to a finite set $C^\prime$ of non-realised cuts in 
$A$ such that $C_A^\theta\nsubseteq C^\prime$ and, for any tuples $\bar{b}$ and 
$\bar{b}'$ from $A$, $otp(\bar{b}/C')=otp(\bar{b}^\prime/C')\implies 
N'\models \theta(\bar{b})\leftrightarrow\theta(\bar{b}^\prime)$. This would be a 
contradiction. 

As each $C_B^\theta$ is finite, it follows that we cannot have an infinite 
sequence of small indiscernible extensions $A\subseteq B_1\subseteq 
B_2\subseteq...\subseteq M'$ such that $| 
C_A^\theta|<|C_{B_1}^\theta|<|C_{B_2}^\theta|<...$ (since otherwise the union 
$\bigcup\limits_{n\in\mathbb{N}}B_n$ would be a small indiscernible sequence for which the conclusion of Fact \ref{CScutfact} is false). So then, for 
each $\theta$, we can find a small extension $B$ of $A$ such that, for any small 
extension $B^\prime$ of $B$, $|C_B^\theta|=|C_{B^\prime}^\theta|$ and in fact 
there is a bijection from $C_B^\theta$ to $C_{B^\prime}^\theta$ which sends each 
cut $c\in C_B^\theta$ to the unique $c^\prime\in C_{B^\prime}^\theta$ which 
refines it.

The process of extending $A$ to such a $B$ could be called ``maximising for 
$\theta$". Enumerate all formulas with parameters in $N$, add one to the 
enumerating indices so that only successor ordinals are used and then maximise 
for each formula in turn, taking unions at limit ordinals. We thereby obtain a 
small indiscernible sequence $A^*\subseteq M'$ extending $A$ such that, for any 
formula $\theta$ with parameters in $N$ and any small indiscernible $B\subseteq 
M'$ extending $A^*$, there is a bijection from $C_{A^*}^\theta$ to $C_B^\theta$ 
which sends each $c\in C_{A^*}^\theta$ to the unique $c^\prime\in C_B^\theta$ 
which refines it.  Note that ``formula $\theta$" really means ``formula $\theta(\bar{x})$ where $\bar{x}$ is a tuple of variables in the sort of $A$".

We would like to consider $\bigcup\limits_\theta C_{A^*}^\theta$ and, for any 
small extension $B$ of $A^*$, the bijection taking each $c\in 
\bigcup\limits_\theta C_{A^*}^\theta$ to the unique $c^\prime\in 
\bigcup\limits_\theta C_B^\theta$ which refines it. However, we cannot be sure 
at this stage that such a bijection exists. The problem is that, for some 
$\theta_1,\theta_2$, we might have $|C_{A^*}^{\theta_1}\cap 
C_{A^*}^{\theta_2}|>|C_B^{\theta_1}\cap C_B^{\theta_2}|$. In other words, some 
cuts might coincide in $A^*$ but not in $B$. For any small $B'$ extending a small $B$ 
extending $A^*$, we must have $|C_{A^*}^{\theta_1}\cap 
C_{A^*}^{\theta_2}|\geq|C_B^{\theta_1}\cap C_B^{\theta_2}|\geq|C_{B'}^{\theta_1}\cap C_{B'}^{\theta_2}|$. (To see this note, in the notation of the first inequality, that every cut in $C_B^{\theta_1}\cap C_B^{\theta_2}$ refines one in $C_{A^*}^{\theta_1}\cap C_{A^*}^{\theta_2}$ and that it would contradict the existence of our bijections if two cuts in $C_B^{\theta_1}\cap C_B^{\theta_2}$ were to refine the same cut in $C_{A^*}^{\theta_1}\cap C_{A^*}^{\theta_2}$.)  An ordinal (in 
this case a finite one) cannot be decreased infinitely many times. So, for each 
pair $\theta_1,\theta_2$, we can extend so that $|C_B^{\theta_1}\cap 
C_B^{\theta_2}|$ is minimised. We can enumerate all such pairs of formulas and extend 
appropriately for each one in turn, taking unions at limit ordinals. This 
results in a small extension $A^{**}\subseteq M'$ of $A^*$ with the following 
property. For any small indiscernible $B\subseteq M'$ extending $A^{**}$ and any 
formula $\theta$ with parameters in $N$, let 
$f_B^\theta:C_{A^{**}}^\theta\rightarrow C_B^\theta$ be the bijection which maps 
each $c\in C_{A^{**}}^\theta$ to the unique $c^\prime\in C_B^\theta$ which 
refines it and define $C_B=\bigcup\limits_\theta C_B^\theta$. Then, for each such 
$B$, the union, over all such $\theta$, of the graphs of the functions 
$f_B^\theta$ is the graph of an order-preserving bijection 
$f_B:C_{A^{**}}\rightarrow C_B$.

To simplify notation, let $C=C_{A^{**}}$. Enumerate the elements of $C$ as 
$(c^\alpha)_{\alpha<\kappa}$. For each small indiscernible $B\subseteq M'$ 
extending $A^{**}$ and each $\beta<\kappa$, let $c_B^\beta=f_B(c^\beta)$. We 
build a chain $(B_\alpha)_{\alpha<\kappa}$ of small indiscernible sequences in 
$M'$ extending $A^{**}$ in the following way. Let $\beta<\kappa$ and suppose we 
have formed $B_\alpha$ for all $\alpha<\beta$. Let 
$B_\beta'=A^{**}\cup\bigcup\limits_{\alpha<\beta} B_\alpha$. Consider the cut 
$c^\beta_{B_\beta'}$. There are two cases. 

\begin{enumerate}
\item Suppose it is possible to extend $B^\prime_\beta$ to a small 
indiscernible sequence $B\subseteq M'$ which has elements realising the cut 
$c^\beta_{B_\beta'}$ to each side of $c^\beta_B$. Choose $u_\beta\in B$ 
realising $c_{B^\prime_\beta}^\beta$ below $c_B^\beta$. Choose $v_\beta\in B$ 
realising $c_{B^\prime_\beta}^\beta$ above $c_B^\beta$. Then let 
$B_\beta=B_\beta^\prime\cup\{u_\beta,v_\beta\}$ with the obvious ordering. 
	
	
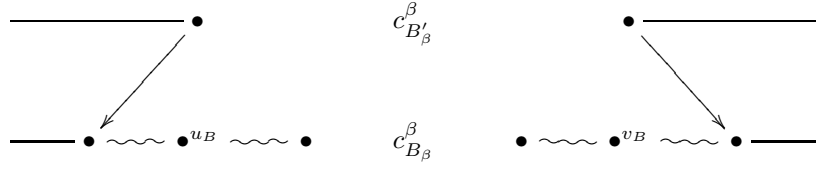
\[
\xymatrix
{
&&\ar@{-}[ll] \bullet\ar[dl]&&c^\beta_{B_\beta'}&&\bullet\ar[dr]\ar@{-}[rr]&& \\
&\ar@{-}[l] \bullet&\bullet^{u_B}\ar@{~}[l]&\bullet 
\ar@{~}[l]&c_{B_\beta}^\beta&\bullet 
\ar@{~}[r]&\bullet^{v_B}\ar@{~}[r]&\bullet\ar@{-}[r]&
}
\]
	
\item If it is not possible to extend $B_\beta'$ as above, choose $B_\beta=B_\beta'$.
\end{enumerate}

Let $A'=\bigcup\limits_{\alpha<\kappa}B_\alpha$. In the remainder of this proof, if we write an indiscernible sequence as $(B,<)$ we are thinking of it as a structure and the language we are using is $\{<\}$. The sequence $A'$ has been constructed so that, whenever we have small indiscernible sequences $B$ and $B'$ in $M'$, with $A'\subseteq B$ and $A'\subseteq B'$, and $\bar{a}=(a_1,...,a_n)\in B^n$ and $\bar{a}'=(a'_1,...,a'_n)\in B'^n$, if $(A',<)\prec (B,<)$ and $(A',<)\prec (B',<)$ and $qftp_{\{<\}}(\bar{a}/A')=qftp_{\{<\}}(\bar{a}'/A')$ then, for each $\beta<\kappa$ and $i\in \{1,...,n\}$, $a_i<c_B^{\beta}$ if and only if $a'_i<c_{B'}^{\beta}$.







We show that $A'$ has the desired property. Let $q(\bar x)$ be a complete $L$-type over $A'$ which is finitely realised in $A'$. 
Let $\bar a=(a_1,...,a_n)$ and $\bar a'=(a_1',...,a_n')$ both realise $q(\bar x)$ in $M'$. Let $(N'',M'',A'')$ be an elementary extension of $(N',M',A')$, in the language $L\cup \{P,Q\}$ where $P$ is a unary predicate for $M'$ and $Q$ is a unary predicate for $A'$, such that $q(\bar{x})$ is realised in $A''$ by $\bar{a}''$. By the downward L\"owenheim-Skolem theorem, there is a small indiscernible (in the sense of $L$) sequence $B''\subseteq A''$ such that $A'\cup\{\bar{a}\}\subseteq B''$ and $(A',<)\prec(B'',<)$. We then have automorphisms $\sigma$ and $\sigma'$ of the $L\cup\{P\}$-structure $(N'',M'')$ which fix $A'$ pointwise and are such that $\sigma(\bar{a}'')=\bar{a}$ and $\sigma'(\bar{a}'')=\bar{a}'$. Let $B=\sigma(B'')$ and $B'=\sigma'(B'')$. Then $B$ and $B'$ are indiscernible (in the sense of $L$), $(A',<)\prec(B,<)$ and $(A',<)\prec(B',<)$. Since $qftp_{\{<\}}(\bar{a}/A')=qftp_{\{<\}}(\bar{a}'/A')$ we have $a_i<c_B^{\beta}$ if and only if $a'_i<c_{B'}^{\beta}$, for each $\beta<\kappa$ and $i\in\{1,...,n\}$. It follows that $tp(\bar{a}/N)=tp(\bar{a}'/N)$.


\begin{rem}\label{more}
One then gets for free that $q(\bar{x})\cup\{P(x_1),...,P(x_n)\}$ implies a complete 
$L$-type over $N A^\prime$. One way to obtain $N$ and $(N^\prime,M^\prime)$ would be 
to take saturated elementary extensions $M^{Sh}\prec \bar{N}$ and 
$(\bar{N},M^{Sh})\prec (\bar{N}^\prime,\bar{M}^\prime)$, in the languages $L(M^{Sh})$ and $L(M^{Sh})\cup\{P\}$ respectively, and then define $N$, 
$M^\prime$ and $N^\prime$ to be the $L$-reducts of $\bar{N}$, $\bar{M}^\prime$ 
and $\bar{N}^\prime$. In this case, given that $q(\bar{x})\cup\{P(x_1),...,P(x_n)\}$ implies a complete $L$-type over $NA'$, $q(\bar{x})$ implies a maximal consistent collection of sets defined by predicates in $L(M^{Sh})$ using parameters in $A'$. Since $Th(M^{Sh})$ has QE (Fact \ref{QE}), $q(\bar{x})$ implies a complete $L(M^{Sh})$-type over $A'$. 
\end{rem}

\section{Theorem}\label{thmsect}

In this section we prove Theorem \ref{main}, after restating it for convenience. 
Our proof borrows a great deal from the proof in \cite{CS} that distality is 
equivalent to having so-called strong honest definitions (see Proposition 19 and Theorem 21 of \cite{CS}). It is essentially just a stretching 
of that argument to a setting provided by Lemma \ref{extending}. 

\begin{thm}\label{distshimpliesdist}
Let $L$ be a language and $T$ a complete first-order $L$-theory. Let $M\models T$ and let $M^{Sh}$ be the Shelah 
expansion of $M$. Suppose $Th(M^{Sh})$ is distal. Then $T$ is distal.
\end{thm}

\proof As in the proof of Lemma \ref{extending}, the default language is $L$. Suppose, for contradiction, that $T$ is not distal. By Lemma \ref{tech} 
there exist a model $K\models T$, an indiscernible sequence $I+\{b\}+J$ in $K$, 
with $I$ and $J$ infinite and $\{b\}$ a singleton, some $a\in K$ and a 
formula $\phi(x,y)$ such that $I+J$ is indiscernible over $a$, $K\models 
\phi(a,d)$ for all $d\in I+J$ and $K\models \neg\phi(a,b)$. 

Let $M^{Sh}\prec \bar{N}$ and $(\bar{N},M^{Sh})\prec (\bar{N}',\bar{M}')$ be 
saturated elementary extensions in the languages $L(M^{Sh})$ and $L(M^{Sh})\cup \{P\}$ respectively. Let $N,M'$ and $N'$ be the reducts of 
$\bar{N},\bar{M}'$ and $\bar{N}'$ to $L$. We may assume $K\prec M^\prime$. 

By Lemma \ref{extending} and Remark \ref{more}, $I+J$ extends to a small 
indiscernible sequence $A'$ in $M'$ with the property that every complete 
type over $A'$ which is finitely realised in $A'$ implies a complete 
$L(M^{Sh})$-type over $A'$. We may assume $a$ and $b$ are such that $A'\cup\{b\}$ is indiscernible, with $b$ positioned just above $I$, and $A'$ is indiscernible over $a$. (To see this, let $r(x,y)$ be the partial type over $A'$ expressing the desired properties of the pair $(a,b)$. Any finite $r'(x,y)\subseteq r(x,y)$ involves only a finite tuple $cd$ from $A'$, where $c$ is bounded above by an element of $I$ and $d$ lies entirely above $I$. By indiscernibility, an automorphism of $M'$ takes $c$ to a tuple in $I$ and $d$ to a tuple in $J$, establishing that $r'(x,y)$ can be realised. One then uses saturation.)

Now consider the structure $(\bar{N}',\bar{M}',A')$ in the language $L(M^{Sh})\cup\{P,Q\}$. Take sufficiently saturated elementary extensions  
\[(\bar{N}',\bar{M}',A')\prec (\bar{N}'',\bar{M}'',A'')\prec 
(\bar{N}''',\bar{M}''',A''')\prec (\bar{N}'''',\bar{M}'''',A'''').\]

In all cases, let the removal of the bar correspond to taking the $L$-reduct. Let 
$p(x)=tp(a/A'')$. Let $q(y)$ be a complete type over $A''$, where $y$ is a 
single variable in the sort of $A'$, such that $q$ is finitely realised in $A'$. 
We show that $p(x)\cup q(y)$ implies a complete type in $xy$ over $\emptyset$.

Let $q'(y)$ be some extension of $q(y)$ to a complete type over $N''''$ which is 
finitely realised in $A'$. Let $(d_i)_{i\in \mathbb{Z}}$ be a Morley sequence for $q'$ over 
$N'$ in $A''$. (Of course, a Morley sequence is usually indexed by $\mathbb{N}$ but, having obtained such a sequence, one can choose another one indexed by $\mathbb{Z}$ with the same EM-types.)
 Let $d^*$ realise $q$ in $M'''$. Let $(d'_i)_{i\in\mathbb{Z}}$ be a 
Morley sequence in $q'$ over $N'''$ in $A''''$. 
We found the following picture helpful.
\begin{center}
\begin{scaletikzpicturetowidth}{\textwidth}
\begin{tikzpicture}
\draw (0,0) -- (15,0);
\draw (0,0) -- (0,4.5);
\draw (0,4.5) -- (15,4.5);
\draw (15,0) -- (15,4.5);
\draw (2,0) -- (2,4.5);
\draw (6,0) -- (6,4.5);
\draw (9,0) -- (9,4.5);
\draw (11.3,0) -- (11.3,4.5);
\draw (0,2.5) -- (15,2.5);
\draw [ultra thick](3,.9) -- (6,.9);
\draw [red] (8.5,1) ellipse (6cm and .4cm);
\node [below left] at (2,2.5) {$M$};
\node [below left] at (6,2.5) {$M'$};
\node [below left] at (9,2.5) {$M''$};
\node [below left] at (11.3,2.5) {$M'''$};
\node [below left] at (15,2.5) {$M''''$};
\node [below left] at (2,4.5) {$N$};
\node [below left] at (6,4.5) {$N'$};
\node [below left] at (9,4.5) {$N''$};
\node [below left] at (11.3,4.5) {$N'''$};
\node [below left] at (15,4.5) {$N''''$};
\node [above left] at (5.5,.8) {$I+J$};
\node [above right] at (6.3,.7) {$... d_{-1},d_0,d_1...$};
\node [above right] at (9.4,1.75) {$\bullet d^*$};
\node [above right] at (11.3,.7) {$... d'_{-1},d'_0,d'_1...$};
\node [above left, red] at (6,1.35) {$A'$};
\node [above right, red] at (8,1.35) {$A''$};
\node [above right, red] at (10.5,1.4) {$A'''$};
\node [above right, red] at (13.4,1.15) {$A''''$};
\end{tikzpicture}

\end{scaletikzpicturetowidth}
\end{center}

The sequence $(d_i)_{i\in\mathbb{Z}}+d^*+(d'_i)_{i\in \mathbb{Z}}$ is $L(M^{Sh})$-indiscernible. This is because every finite subsequence has a type over $A'$ which is finitely realised in $A'$ and therefore implies a complete $L(M^{Sh})$-type over $A'$. Also $(d_i)_{i\in \mathbb{Z}}+(d'_i)_{i\in\mathbb{Z}}$ is $L$-indiscernible over $Na$ and so 
$L(M^{Sh})$-indiscernible over $a$. By the distality of $Th(M^{Sh})$,  
 $(d_i)_{i\in\mathbb{Z}}+d^*+(d'_i)_{i\in \mathbb{Z}}$ is 
$L(M^{Sh})$-indiscernible over $a$ and thus $L$-indiscernible over $a$. Since 
$d^*$ was an arbitrary realisation of $q$ in $M'''$ it follows that 
$p(x)\cup q(y)$ implies a complete type in $xy$ over $\emptyset$. 

Since the set of all $d \in M'''$ such that $tp(d/A'')$ is finitely realised in $A'$ is type-definable over $A''$ in the structure $M'''$, a
compactness argument gives us some $\bar{c}\in A''^k$ and an $L$-formula $\theta(x, \bar{z})$, 
with $\bar{z}$ a $k$-tuple of variables in the sort of $y$, such that $M''\models \theta (a,\bar{c})$ and 
$\theta (x,\bar{c})$ implies the $\phi$-type of $a$ over $A'$. (To see this, suppose not. Then, for any choice of $\theta(x,\bar{c})\in tp(a/A'')$, there exist $b'\in A'\subseteq\{d\in M'':tp(d/A'')\text{ is finitely realised in }A'\}$ and $a'\in M''$ such that $M''\models \theta(a',\bar{c})$ and $M''\models\phi(a,b')\leftrightarrow\neg\phi(a',b')$. Then, by compactness, there exist $a',b'\in M''$ such that $a'\models p(x)$, $tp(b'/A'')$ is finitely realised in $A'$ and $tp(ab')\neq tp(a'b')$ which is a contradiction.) So for any finite 
$A\subseteq A'$ there is a $\bar{c}\in A'^k$ such that $M'\models \theta (a,\bar{c})$ and 
$\theta (x, \bar{c})$ implies the $\phi$-type of $a$ over $A$. 

Recall that $A'\cup \{b\}$ is indiscernible, with $b$ positioned just above $I$, and $A'$ is 
indiscernible over $a$. So $M'\models \phi (a,d)$ for all $d\in A'$ and  
$M'\models \neg\phi (a,b)$. Let $A\subseteq A'$ have cardinality $k+1$. Let 
$\bar{c}\in A'^k$ be such that $M'\models \theta (a, \bar{c})$ and $\theta (x, \bar{c})$ implies 
the $\phi$-type of $a$ over $A$. Let $d$ be an element of $A$ which does not 
belong to the tuple $\bar{c}=(c_1,...,c_k)$. Let $f$ be a partial automorphism of $(A'\cup\{b\},<)$, 
with domain $A \cup\{c_1,...,c_k\}$, such that $f(d)=b$. Let $f(\bar{c})$ denote the tuple $(f(c_1),...,f(c_k))$. Then $M'\models \theta(a,f(\bar{c}))$. It follows that $M'\models \phi(a,b)$, since otherwise there would be some $a'\in M'$ such that $M'\models \theta(a',\bar{c})\wedge\neg\phi(a',d)$ which would be a contradiction. But $M'\models\neg\phi(a,b)$ and so we have a contradiction and the proof is finished. \endproof

\section{Converse}

For completeness we give a proof of the converse of Theorem \ref{main}, though 
it was already known to experts in the area. 

\begin{thm}\label{distimpliesdistsh}
Let $T$ be a complete first-order $L$-theory and let $M\models T$. Suppose $T$ is distal. Then $Th(M^{Sh})$ is distal. 
\end{thm}

\proof Let $M^{Sh}\prec\bar{N}$ and $(\bar{N},M^{Sh})\prec(\bar{N}',\bar{M}')$ 
be sufficiently saturated elementary extensions, the first in the language $L(M^{Sh})$ and the second in the language $L(M^{Sh})\cup\{P\}$, where $P$ is a new unary predicate. Let $N,M'$ and $N'$ be the 
$L$-reducts of $\bar{N},\bar{M}'$ and $\bar{N}'$ respectively. Let $I+\{b\}+J$ 
be a small $L(M^{Sh})$-indiscernible sequence in $\bar{M}'$, such that $I$ and $J$ are both infinite without endpoints. Let $A\subseteq \bar{M}'$ be 
small and suppose $I+J$ is $L(M^{Sh})$-indiscernible over $A$. It follows that 
$I+J$ is $L$-indiscernible over $NA$ and that $I+\{b\}+J$ is $L$-indiscernible over 
$N$. By distality of $T$, $I+\{b\}+J$ is $L$-indiscernible over $NA$. Therefore 
$I+\{b\}+J$ is $L(M^{Sh})$-indiscernible over $A$. \endproof

\section{Alternative approach using measures}

In this section we mention an alternative proof of Theorem \ref{main} for which we thank Ehud Hrushovski and Anand Pillay. We shall be brief with the details as we have already given a thorough proof and experts will be able to fill in the gaps quite easily. Throughout, $T$ is a complete first-order NIP theory. 

We recall the definitions we shall need and direct the reader to Chapter 7 of \cite{Sbook} for further details. A measure $\mu$ over a structure $M$ assigns to each definable (with parameters) subset of some fixed sort of $M$ a number in the interval $[0,1]$. It is finitely additive and achieves a maximum value of $1$. It is smooth if, for each elementary extension $M\prec N$, there is only one measure over $N$ extending $\mu$. Suppose we have $M\models T$, a sufficiently saturated elementary extension $M\prec N$ and a measure $\mu$ over $N$. We say $\mu$ is definable over $M$ if, for each $L$-formula $\varphi(x,y)$ and closed $B\subseteq [0,1]$, the set of all $b\in N$ such that $\varphi(x,b)$ defines a set with $\mu$-measure in $B$ is type-definable over $M$. We say $\mu$ is finitely satisfiable in $M$ if every definable set with positive $\mu$-measure has non-empty intersection with the relevant sort of $M$. In the case where $\mu$ is both definable over $M$ and finitely satisfiable in $M$, we say $\mu$ is generically stable with respect to $M$.

The following result goes back to \cite{Sdist} but perhaps it is most convenient to refer the reader to Propositions 9.26 and 9.27 of \cite{Sbook}. 

\begin{fact}\label{Pierre}
Let $T$ be a complete first-order $L$-theory. Suppose $T$ has NIP. Then $T$ is distal if and only if, for every $M\models T$, sufficiently saturated $M\prec N$ and measure $\mu$ over $N$, if $\mu$ is generically stable with respect to $M$ then its restriction to $M$ is smooth. 
\end{fact}

We can obtain Theorem \ref{main} by combining this with the following.  

\begin{thm}\label{UdiAnand}
Let $M\models T$. Let $M^{Sh}\prec \bar K$. Then let $\bar K\prec \bar{N}$ and $(\bar{N},\bar K)\prec (\bar{N}',\bar{K}')$ be sufficiently saturated elementary extensions and define $N,K, N'$ and $K'$ to be the $L$-reducts of $\bar{N},\bar{K}, \bar {N}'$ and $\bar{K}'$. Let $\mu$ be a measure over $N'$ which is generically stable with respect to $K$. Then the restriction of $\mu$ to $K$ extends to a measure $\mu^*$ over $\bar K$ with the following properties:

\begin{enumerate}

\item $\mu^*$ extends to a measure over $\bar{K}'$ which is generically stable with respect to $\bar K$ and

\item if the restriction of $\mu$ to $K$ is not smooth then $\mu^*$ is not smooth.

\end{enumerate}

Furthermore, the measure $\mu^*$ is the unique extension to $\bar K$ of the restriction of $\mu$ to $K$. 

\end{thm}

\proof Let $X$ be a definable set of $\bar{K}'$. Then $X$ is a fibre of a $\emptyset$-definable set, say $X_1$. By Fact \ref{QE}, $X_1$ is defined by a formula $R(xz)$ where $R\in L(M^{Sh})$. Consider the set, say $X_2$, defined by $R(xz)$ in $K$. As is well known, $X_2$ will itself be externally definable with respect to the structure $K$ and so there will exist an $L$-formula $\varphi(x,y,z)$ and parameter $b$ from $N$ such that $X_2$ is the set of all $(a,c)\in K$ for which $N\models\varphi(a,b,c)$. Then $X_1$ will be the set externally defined by $\varphi(x,b,z)$ in $K'$. So, finally, we have an $L$-formula $\varphi(x,y,z)$ and parameters $b$ from $N$ and $c$ from $K'$ such that $X=\{a\in K':N'\models\varphi(x,b,c)\}$. Define $\mu'$ over $\bar{K}'$ such that $\mu'(X)$ is the value assigned by $\mu$ to the set defined by $\varphi(x,b,c)$ in $N'$. 

One checks that $\mu'$ is a measure over $\bar{K}'$. Note that it is well-defined because $\varphi(x,b,c)$ is always unique up to a $\mu$-measure zero symmetric difference, using the fact that $\mu$ is finitely realisable in $K$. Furthermore, $\mu'$ is definable over $\bar K$ and finitely satisfiable in $\bar K$. For definability, we use definability of $\mu$ for the formula $\varphi(x,y,z)$ and then restrict the type-definable set to $y=b$. Finite realisability of $\mu'$ in $\bar{K}$ is immediate from the finite realisability of $\mu$ in $K$. So then the measure $\mu'$ over $\bar{K}'$ is generically stable with respect to $\bar K$. We define $\mu^*$ to be its restriction to $\bar K$.  

For uniqueness of $\mu^*$ (the ``furthermore" statement) one uses the fact that every measure over $\bar K$ comes from a measure over $N$ which agrees with it on $K$ and is finitely realisable in $K$. Since $\mu$ is generically stable with respect to $K$, it is known (by Proposition 3.3 in \cite{NIP3}) to have only one finitely realisable (in $K$) extension to $N$. Therefore $\mu^*$ is unique. 

Now suppose the restriction of $\mu$ to $K$ is not smooth. Since $K\prec N$ is sufficiently saturated, this is witnessed over $N$ and so there are two distinct extensions $\mu_1$ and $\mu_2$ to $N$. These extend, respectively, to measures $\mu_1^*$ and $\mu_2^*$ over $\bar{N}$. Trivially, $\mu_1^*$ and $\mu_2^*$ are distinct. They both restrict to measures on $\bar K$ which extend the restriction of $\mu$ to $K$. By uniqueness of $\mu^*$, they both extend $\mu^*$. So $\mu^*$ is not smooth. \endproof

One quickly deduces Theorem \ref{main} from this as follows. Suppose $T$ is not distal and let $M\models T$. Let $L_M$ be $L$ together with a new constant symbol for every element of $M$. Let $T_M$ be the $L_M$-theory of $M$. By Corollary 2.9 in \cite{Sdist}, $T_M$ is not distal. Then, by Fact \ref{Pierre}, there exist a model $\hat K\models T_M$, a sufficiently saturated elementary extension $\hat K\prec \hat N'$ and a measure $\hat \mu$ over $\hat N'$ which is generically stable with respect to $\hat K$ and whose restriction to $\hat K$ is not smooth. Let $N'$ and $K$ be the $L$-reducts of $\hat N'$ and $\hat K$ and let $\mu$ be the restriction of $\hat \mu$ to $N'$. Then $\mu$ is generically stable over $K$ and its restriction to $K$ is not smooth. We may assume $N'$, $K$ and $\mu$ are as in Theorem \ref{UdiAnand} and so there exist $\bar K, \bar{K}'$ and $\bar{N}'$ as in Theorem \ref{UdiAnand}. It follows that there is a non-smooth measure $\mu^*$ over $\bar{K}$ which extends to a measure over $\bar K'$ which is generically stable with respect to $\bar K$. By Fact \ref{Pierre}, $Th(\bar K)$ is not distal. So $Th(M^{Sh})$ is not distal.

\bibliographystyle{plain}

\bibliography{referencesBKup}

\begin{thebibliography}{1}

\bibitem{CS}
Artem Chernikov and Pierre Simon.
\newblock Externally definable sets and dependent pairs ii.
\newblock {\em Transactions of the American Mathematical Society}, 367:5217 --
  5235, 2015.

\bibitem{NIP3}
Ehud Hrushovski, Anand Pillay, and Pierre Simon.
\newblock Generically stable and smooth measures in {NIP} theories.
\newblock {\em Trans. Amer. Math. Soc.}, 365(5):2341--2366, 2013.

\bibitem{ShelahMB}
Saharon Shelah.
\newblock Classification theory for elementary classes with the dependence
  property---a modest beginning.
\newblock {\em Sci. Math. Jpn.}, 59(2):265--316, 2004.
\newblock Special issue on set theory and algebraic model theory.

\bibitem{ShExp}
Saharon Shelah.
\newblock Dependent first-order theories, continued.
\newblock {\em Israel Journal of Mathematics}, 173:1 -- 60, 2009.

\bibitem{NIPtypedecomp}
Pierre Simon.
\newblock Type decomposition in nip theories.
\newblock To appear in {\it Journal of the European Mathematical Society}.

\bibitem{Sdist}
Pierre Simon.
\newblock Distal and non-distal {NIP} theories.
\newblock {\em Ann. Pure Appl. Logic}, 164(3):294--318, 2013.

\bibitem{Sbook}
Pierre Simon.
\newblock {\em A Guide to NIP Theories}.
\newblock Cambridge University Press, 2015.
\newblock Cambridge Books Online.

\end{thebibliography}

\end{document}